\documentclass[11pt, reqno]{amsart}
\usepackage{preamble}
\newcommand{\T}{\mathbb{T}}
\newcommand{\cD}{\mathcal{D}}
\newcommand{\Var}{\mathrm{Var}}
\newcommand{\Cov}{\mathrm{Cov}}
\allowdisplaybreaks

\title{Improved bounds for non-vanishing transversals}
\author{Adrian Beker}
\address{Department of Mathematics, Stanford University, Stanford, CA 94305}
\email{abeker@stanford.edu}

\begin{document}

\begin{abstract}
    Given $A_1,\ldots,A_m \subseteq \Z$, a \emph{non-vanishing transversal} is a sequence $(a_1,\ldots,a_m) \in A_1\times\ldots\times A_m$ with no non-empty zero-sum subsequence. Let $g(m)$ be the least positive integer $n$ such that any $A_1,\ldots,A_m$ of size $n$ have a non-vanishing transversal. It is not difficult to see that $g(m) \geq m+1$, and Goddyn and Tarsi conjectured that this is tight. In this direction, Alon and Ruzsa proved the upper bound $g(m) \leq m^{1+o(1)}$, in the quantitative form $g(m) \leq m\exp\bigl(O(\sqrt{\log m\log\log m})\bigr)$. In this paper, we improve this to $g(m) = O\bigl(m(\log m)^2\bigr)$. The key ingredient is an improved upper bound on the expected minimum fractional part in a random dilate of a set of $n$ integers, which was found by GPT-6 Astra.
\end{abstract}

\maketitle

\section{Introduction}\label{sec:intro}

In this paper, we are interested in the following problem in combinatorial number theory. Given sets $A_1,\ldots,A_m \subseteq \Z$, a \emph{transversal} is a sequence $(a_1,\ldots,a_m)$ such that $a_i \in A_i$ for each $i \in [m]$. Such a transversal is said to be \emph{non-vanishing} if $\sum_{i \in I}a_i \neq 0$ for all non-empty $I \subseteq [m]$. Given $m$, how large does the smallest of the sets $A_1,\ldots,A_m$ need to be in order to guarantee the existence of a non-vanishing transversal? Letting $g(m)$ denote the answer to this question, it is not difficult to see that $g(m) > m$. Indeed, let $B = \{b_1,\ldots,b_m\}$ be any set of $m$ integers, and for $i \in [m]$ define
\[
    A_i \vcentcolon= B - b_i = \{b_j - b_i \mid j \in [m]\}.
\]
Any transversal has the form $a_i = b_{f(i)} - b_i$ for some function $f \colon [m] \to [m]$. But then a cycle in the functional graph of $f$ gives rise to a zero-sum subsequence.

In the other direction, a simple greedy algorithm shows that $g(m) \leq 2^{m-1}+1$. Indeed, if $|A_i| > 2^{m-1}$ for all $i \in [m]$, we can construct the desired transversal by iteratively choosing $a_i$ to be an element of $A_i$ not of the form $-\sum_{j\in J}a_j$ for $J \subseteq [i-1]$. However, this upper bound is far from tight: Alon and Ruzsa \cite{alon-ruzsa} showed that $g$ is in fact close to linear, in the sense that $g(m) \leq m^{1+o(1)}$. More precisely, they proved the following result.

\begin{theorem}
    \label{thm:main_old}
    $g(m) \leq me^{O(\sqrt{\log m\log\log m})}$
\end{theorem}

Goddyn and Tarsi \cite[Conjecture 1.2]{alon-ruzsa} conjectured that $g(m) = m+1$, i.e.\ that the above example is in fact best possible. This problem was also mentioned in the collection \emph{Open Problems in Additive Combinatorics} compiled by Croot and Lev \cite[Problem 6.1]{croot-lev}. 

Although we are not able to resolve the conjecture of Goddyn and Tarsi, we make further progress by establishing the following quantitative improvement over Theorem \ref{thm:main_old}.

\begin{theorem}
    \label{thm:main_new}
    $g(m) = O\bigl(m(\log m)^2\bigr)$
\end{theorem}

In their proof of Theorem \ref{thm:main_old}, Alon and Ruzsa use a probabilistic argument to reduce the problem to finding an upper bound for the expected minimum fractional part in a random dilate of a set of $n$ integers. We establish Theorem \ref{thm:main_new} by making the following improvement to this bound, which may be of independent interest. The proof of this result was found by GPT-6 Astra (see Section \ref{sec:ai} for more details).

\begin{theorem}
    \label{thm:main_ingredient}
    Let $A$ be a set of $n$ non-zero integers. Then\footnote{Here and throughout, $\{x\}$ denotes the fractional part of $x \in \R$.}
    \begin{equation*}
        \int_{0}^{1}\min_{a\in A}\{a\theta\}\,d\theta \ll \frac{(\log n)^2}{n}.
    \end{equation*}
\end{theorem}

In the rest of the paper, we will assume that the sets $A_1,\ldots,A_m$ in the definition of $g(m)$ do not contain zero. This changes the value of $g(m)$ by exactly $1$, since zero cannot be part of a non-vanishing transversal anyway. As we are only interested in the asymptotic behavior of $g$, this modification will not affect any of the results.

The paper is organized as follows. In Section \ref{sec:outline}, we first recall the approach of Alon and Ruzsa in more detail. We then give an outline of our approach, compare it with that of Alon and Ruzsa, and highlight the key differences. Sections \ref{sec:int_cov}-\ref{sec:pf_of_thm} contain the proofs of the intermediate results referred to in Section \ref{sec:outline}. In Section \ref{sec:conc_rem}, we make some concluding remarks and comment on the limitations of our approach.

\subsection{Notation}

Our notational conventions are largely standard. As usual, for a positive integer $n$, we abbreviate the set $\{1,\ldots,n\}$ to $[n]$. We use Vinogradov asymptotic notation, so given quantities $A$ and $B$, we write $A \ll B$ to mean $|A| \leq C|B|$ for some absolute constant $C > 0$; this is equivalent to the notation $B \gg A$. If $A \ll B$ and $A \gg B$ hold simultaneously, we write $A \asymp B$. If one of $A$, $B$ is a function on the positive integers, we implicitly assume the underlying argument is sufficiently large. Logarithms are to base $e$ unless otherwise stated.

\subsection{Acknowledgements}

The author would like to thank Jacob Fox and Carl Schildkraut for helpful remarks.

\section{Outline of the argument}\label{sec:outline}

\subsection{The Alon--Ruzsa approach}

The starting point of the approach of Alon and Ruzsa \cite{alon-ruzsa} is to transfer the problem to the torus $\T = \R/\Z$. Indeed, identifying $\T$ with the interval $[0,1)$, note that if $x_1,\ldots,x_m$ are elements of $\T$ such that $x_i > 0$ for all $i \in [m]$ and $\sum_{i=1}^{m}x_i < 1$, then $\sum_{i\in I}x_i$ cannot be zero for any non-empty $I \subseteq [m]$. The idea then is to consider the images of the sets $A_1,\ldots,A_m$ under a suitably chosen homomorphism $\Z \to \T$. More specifically, it suffices to find $\theta \in [0,1)$ such that $\sum_{i=1}^{m}f_{A_i}(\theta) < 1$ and all summands are non-zero, where we define
\[
    f_A(\theta) \vcentcolon= \min_{a\in A}\{a\theta\}.
\]
To this end, Alon and Ruzsa make a random choice of $\theta$. Since $f_{A_i}(\theta) = 0$ happens with probability zero for each $i \in [m]$, it is enough to ensure that $\sum_{i=1}^{m}f_{A_i}(\theta) < 1$ with positive probability. To establish this, it in turn suffices to show that the expected value of $\sum_{i=1}^{m}f_{A_i}(\theta)$ is less than $1$. By linearity of expectation, this amounts to showing that the expected value of each summand is less than $1/m$. For this purpose, Alon and Ruzsa prove the following result; inverting the dependence between $m$ and $n$ then yields Theorem \ref{thm:main_old}. Here and throughout, $U$ denotes a random variable with uniform distribution in $[0,1)$.

\begin{theorem}[{\cite[Lemma 3.2]{alon-ruzsa}}]
    \label{thm:min_fract_part_old}
    Let $A$ be a set of $n$ non-zero integers. Then
    \[
        \E\bigl[f_A(U)\bigr] \leq \frac{e^{O(\sqrt{\log n\log\log n})}}{n}.
    \]
\end{theorem}

Towards a proof of Theorem \ref{thm:min_fract_part_old}, observe that the expectation on the left-hand side can be expressed using layer cake representation as
\begin{equation}
    \label{eq:layer_cake_first}
    \int_{0}^{1}\P\bigl(f_A(U)\geq t\bigr)\,dt = \int_{0}^{1}\P\Bigl(\bigcap_{a\in A}\bigl\{\{aU\} \geq t\bigr\}\Bigr)\,dt.
\end{equation}
At this point it is convenient to assume that $A \subseteq \N$, and we may do so without significant losses. Indeed, at least $n/2$ elements of $A$ have the same sign, and for $a \in \N$ we have $\{-aU\} = \{aU'\}$ where $U' \vcentcolon= 1-U$ is again a uniform random variable in $[0,1)$. Hence, defining for $a \in \N$ and $t \in [0,1)$ the set
\begin{equation}
    \label{eq:interval_union}
    E_{a,t} \vcentcolon= \bigcup_{j=0}^{a-1}\Bigl[\frac{j}{a}, \frac{j+t}{a}\Bigr),
\end{equation}
the event that $\{aU\} < t$ corresponds precisely to $U$ lying in the set $E_{a,t}$. Thus, letting $\lambda$ denote Lebesgue measure on $[0,1)$, it follows that
\begin{equation}
    \label{eq:layer_cake_second}
    \E\bigl[f_A(U)\bigr] = \int_{0}^{1}\lambda\Bigl(\bigcap_{a\in A}E_{a,t}^c\Bigr)\,dt.
\end{equation}
In order to bound the above integrand, Alon and Ruzsa interpret it as the probability that the counting random variable 
\begin{equation}
    \label{eq:counting}
    Z_{A,t} \vcentcolon= \sum_{a\in A}1_{E_{a,t}}(U)
\end{equation}
equals zero, and apply the second moment method. Specifically, Chebyshev's inequality yields
\begin{equation}
    \label{eq:chebyshev}
    \P(Z_{A,t} = 0) \leq \frac{\Var(Z_{A,t})}{\E[Z_{A,t}]^2}.
\end{equation}
Crucially, the first and second moment of $Z_{A,t}$ are tractable, and one can obtain rather precise estimates for them. Indeed, by linearity of expectation,
\begin{equation}
    \label{eq:expectation}
    \E[Z_{A,t}] = \sum_{a\in A}\P\bigl(U \in E_{a,t}\bigr) = \sum_{a\in A}\lambda(E_{a,t}) = nt.
\end{equation}
Furthermore, by bilinearity of covariance,
\begin{equation}
    \label{eq:variance}
    \Var(Z_{A,t}) = \sum_{a_1,a_2 \in A}\Cov\bigl(1_{E_{a_1,t}}(U),1_{E_{a_2,t}}(U)\bigr).
\end{equation}
To deal with the covariance of $E_{a_1,t}$ and $E_{a_2,t}$, Alon and Ruzsa prove the following lemma, which allows for some extra flexibility in the choice of parameters. Here and in the rest of the paper, for $x,y \in [0,1)$ we write
\[
    \Delta(x,y) \vcentcolon= \min(x,y) - xy.
\]
Note that this quantity can be interpreted as the covariance between the events $\{U < x\}$ and $\{U < y\}$, and is always non-negative.

\begin{lemma}[{\cite[Lemma 3.5]{alon-ruzsa}}]
    \label{lm:covariance}
    Let $a_1,a_2 \in \N$, $t_1,t_2 \in [0,1)$ and write $a_i' \vcentcolon = a_i/\gcd(a_1,a_2)$ for $i \in [2]$. Then
    \[
        \Cov\bigl(1_{E_{a_1,t_1}}(U),1_{E_{a_2,t_2}}(U)\bigr)  = \lambda\bigl(E_{a_1,t_1}\cap E_{a_2,t_2}\bigr) - t_1t_2 = \frac{\Delta\bigl(\{a_1't_2\}, \{a_2't_1\}\bigr)}{a_1'a_2'}.
    \]
\end{lemma}

Thus, combining \eqref{eq:chebyshev}, \eqref{eq:expectation}, \eqref{eq:variance} and using Lemma \ref{lm:covariance} to estimate
\begin{equation}
    \label{eq:covariance}
    \Cov\bigl(1_{E_{a_1,t}}(U),1_{E_{a_2,t}}(U)\bigr) \leq \frac{t\min(a_1',a_2')}{a_1'a_2'} = \frac{t\gcd(a_1,a_2)}{\max(a_1,a_2)},
\end{equation}
we obtain that
\[
    \P(Z_{A,t} = 0) \leq \frac{1}{n^2t}\sum_{a_1,a_2\in A}\frac{\gcd(a_1,a_2)}{\max(a_1, a_2)}.
\]
The final and most substantial task is to estimate the resulting arithmetic sum. Alon and Ruzsa accomplish this using an ingenious combination of graph-theoretic and elementary number-theoretic techniques. In this way, they establish the following.

\begin{theorem}[{\cite[Lemma 3.10]{alon-ruzsa}}]
    \label{thm:gal_sum_alon_ruzsa}
    Let $A$ be any set of $n$ non-zero integers. Then
    \[
        \sum_{a_1,a_2\in A}\frac{\gcd(|a_1|,|a_2|)}{\max(|a_1|, |a_2|)} \leq ne^{O(\sqrt{\log n\log\log n})}.
    \]
\end{theorem}

Going back to the integrand in \eqref{eq:layer_cake_second} and employing the estimate coming from Theorem \ref{thm:gal_sum_alon_ruzsa} when $t \geq 1/n$, and the trivial estimate otherwise, one obtains Theorem \ref{thm:min_fract_part_old}.

Alon and Ruzsa remark that Theorem \ref{thm:gal_sum_alon_ruzsa} is nearly tight and provide an example of a set $A$ for which the left-hand side is as large as $n\exp\bigl(\Omega(\sqrt{\log n/\log\log n})\bigr)$. As an aside, it is worth noting that this sum is closely related to an instance of so-called Gál sums, also known as GCD sums. These sums arise in the study of large values of the Riemann zeta function on the critical line, and in the case of interest, sharp bounds have recently become available. Specifically, de la Bretèche and Tenenbaum \cite{delabreteche-tenenbaum} proved that
\begin{equation}
    \label{eq:gal_sums}
    \sup_{A \subseteq \N,\ |A| = n}\sum_{a_1,a_2\in A}\frac{\gcd(a_1,a_2)}{\sqrt{a_1a_2}} = n\exp\Biggl(\bigl(2\sqrt{2}+o(1)\bigr)\sqrt{\frac{\log n\log\log\log n}{\log\log n}}\Biggr).
\end{equation}
However, incorporating this into the argument of Alon and Ruzsa gives only a marginal improvement over Theorem \ref{thm:main_old}, and new ideas are necessary to push the bound further.

\subsection{The new approach}

Our approach broadly follows that of Alon and Ruzsa, but we make several crucial points of departure. The principal new ingredient is Theorem \ref{thm:main_ingredient}, which improves on Theorem \ref{thm:min_fract_part_old}. We repeat its statement here for the reader's convenience.

\begin{theorem}
    \label{thm:min_fract_part_new}
    Let $A$ be a set of $n$ non-zero integers. Then
    \begin{equation}
        \label{eq:min_frac_part_new}
        \E\bigl[f_A(U)\bigr] \ll \frac{(\log n)^2}{n}.
    \end{equation}
\end{theorem}

Theorem \ref{thm:main_new} follows from Theorem \ref{thm:min_fract_part_new} in the same way in which Theorem \ref{thm:main_old} was deduced from Theorem \ref{thm:min_fract_part_old} in the previous subsection.

The main inefficiency in the previous argument is the substantial correlation between the intervals constituting the sets $E_{a,t}$ for different $a \in A$. Specifically, suppose a rational $q \in [0,1)$ is the left endpoint of one or more intervals $\bigl[j/a,(j+t)/a\bigr)$ with $a \in A$ and $j \in \{0,\ldots,a-1\}$. The overlap between these intervals is redundant as regards the union $\bigcup_{a\in A}E_{a,t}$, but contributes significantly towards the aforementioned correlation. To get around this issue, we retain only the longest such interval and discard the rest. To implement this modification, note that the possible denominators of such a rational $q$ in lowest terms are exactly the divisors of the elements of $A$. For any such divisor $d$, only the least multiple of $d$ in $A$ is relevant, since it yields the longest intervals starting at reduced fractions with denominator $d$.

The above considerations motivate us to define, for a set $A \subseteq \N$,
\[
    D(A) \vcentcolon= \bigl\{d \in \N : \exists a \in A \text{ such that } d \mid a\bigr\}, \quad m_A(d) \vcentcolon= \min\bigl\{a \in A : d \mid a\bigr\},
\]
the latter being valid for any $d \in D(A)$. We also define, for $d \in D(A)$ and $t \in [0,1)$, the set
\[
    F_{d,t}^A \vcentcolon= \bigcup_{\substack{0 \leq j < d\\\gcd(j,d)=1}}\Bigl[\frac{j}{d},\frac{j}{d}+\frac{t}{m_A(d)}\Bigr).
\]
By the previous discussion, we have the key property
\[
    \bigcup_{a\in A}E_{a,t} = \bigcup_{d\in D(A)}F_{d,t}^A.  
\]
Similarly as for the left-hand side, we may express the measure of the complement of the right-hand side as $\P(W_{A,t} = 0)$, where now
\begin{equation}
    \label{eq:counting_coprime}
    W_{A,t} \vcentcolon= \sum_{d\in D(A)}1_{F_{d,t}^A}(U).
\end{equation}
To proceed with the second moment method as before, we need information about the mean and variance of $W_{A,t}$. The former can be easily obtained via linearity of expectation:
\[
    \E[W_{A,t}] = \sum_{d\in D(A)}\P\bigl(U \in F_{d,t}^A\bigr) = \sum_{d\in D(A)}\lambda(F_{d,t}^A) = H(A)t,
\]
where we define
\[
    H(A) \vcentcolon= \sum_{d\in D(A)}\frac{\varphi(d)}{m_A(d)}.
\]
Here and throughout, $\varphi$ denotes Euler's totient function.

To calculate the variance of $W_{A,t}$, we take a step back and relate $W_{A,t}$ to a variant of the random variable $Z_{A,t}$ defined in \eqref{eq:counting}. Specifically, defining for $a \in \N$ and $\alpha \geq 0$ the functions
\begin{equation}
    \label{eq:chi}
    \chi_{a,\alpha} \vcentcolon= \sum_{j \in \Z}1_{[j/a,j/a+\alpha)}, \quad \chi_{a,\alpha}^* \vcentcolon= \sum_{\substack{j \in \Z\\\gcd(j,a)=1}}1_{[j/a,j/a+\alpha)},
\end{equation}
a pointwise application of the Möbius inversion formula yields
\[
    \chi_{a,\alpha}^* = \sum_{d\mid a}\mu(a/d)\chi_{d,\alpha}.
\]
Noting that $1_{F_{d,t}^A}(U) = \chi_{d,t/m_A(d)}^*(U)$ and exploiting the fact that covariance is bilinear, it follows that
\begin{align}
    \Var(W_{A,t}) &= \sum_{d_1,d_2 \in D(A)}\Cov\bigl(1_{F_{d_1,t}^A}(U), 1_{F_{d_2,t}^A}(U)\bigr)\nonumber\\
    \label{eq:var_expansion}
    &= \sum_{d_1,d_2\in D(A)}\sum_{\substack{e_1\mid d_1\\e_2\mid d_2}}\mu\Bigl(\frac{d_1}{e_1}\Bigr)\mu\Bigl(\frac{d_2}{e_2}\Bigr)\Cov\bigl(\chi_{e_1,t/m_A(d_1)}(U), \chi_{e_2,t/m_A(d_2)}(U)\bigr).
\end{align}
Since $\chi_{a,\alpha}(U) = 1_{E_{a,\alpha a}}(U)$ for $\alpha < 1/a$, applying Lemma \ref{lm:covariance} yields the expression
\begin{equation}
    \label{eq:var_identity}
    \Var(W_{A,t}) = \sum_{d_1,d_2\in D(A)}\sum_{\substack{e_1\mid d_1\\e_2\mid d_2}}\mu\Bigl(\frac{d_1}{e_1}\Bigr)\mu\Bigl(\frac{d_2}{e_2}\Bigr)\frac{\gcd(e_1,e_2)^2}{e_1e_2}\Delta\Bigl(\Bigl\{\frac{\lcm(e_1,e_2)t}{m_A(d_1)}\Bigr\}, \Bigl\{\frac{\lcm(e_1,e_2)t}{m_A(d_2)}\Bigr\}\Bigr).
\end{equation}
The above sum typically contains negative terms, so unlike in the Alon--Ruzsa argument, we cannot proceed by merely upper bounding $\Delta$. Instead, we integrate out the covariances before estimating them. To this end, it is convenient to extend the range of integration to infinity. By appropriately extending the definition of $W_{A,t}$ to the range $t \geq 1$, the above formula for the variance remains valid, and one obtains
\[
    \E[f_A(U)] = \int_{0}^{1}\P(W_{A,t} = 0)\,dt \leq \varepsilon + \int_{\varepsilon}^{\infty}\frac{\Var(W_{A,t})}{\E[{W_{A,t}}]^2}\,dt,
\]
where $\varepsilon > 0$ is a (small) suitable threshold. By the above calculations, the last integral can be expressed as
\begin{equation}
    \label{eq:int_var}
    \frac{1}{H(A)^2}\sum_{d_1,d_2\in D(A)}\sum_{\substack{e_1\mid d_1\\e_2\mid d_2}}\mu\Bigl(\frac{d_1}{e_1}\Bigr)\mu\Bigl(\frac{d_2}{e_2}\Bigr)\frac{\gcd(e_1,e_2)^2}{e_1e_2}\int_{\varepsilon}^{\infty}\Delta\Bigl(\Bigl\{\frac{\lcm(e_1,e_2)t}{m_A(d_1)}\Bigr\}, \Bigl\{\frac{\lcm(e_1,e_2)t}{m_A(d_2)}\Bigr\}\Bigr)\Big/t^2\,dt.
\end{equation}
By a fairly routine calculation involving a polarization-type identity for $\Delta$ and Stirling's approximation, the above integral can be evaluated explicitly; we refer the reader to Section \ref{sec:int_cov} for the details. The upshot is that one ends up with a sum over $d_1,d_2 \in D(A)$ weighted by the following two arithmetic quantities:
\begin{align}
    \label{eq:arithm_sums_g}
    G(d_1,d_2) &\vcentcolon= \sum_{\substack{e_1\mid d_1\\e_2\mid d_2}}\mu\Bigl(\frac{d_1}{e_1}\Bigr)\mu\Bigl(\frac{d_2}{e_2}\Bigr)\gcd(e_1,e_2),\\
    \label{eq:arithm_sums_l}
    L(d_1,d_2) &\vcentcolon= \sum_{\substack{e_1\mid d_1\\e_2\mid d_2}}\mu\Bigl(\frac{d_1}{e_1}\Bigr)\mu\Bigl(\frac{d_2}{e_2}\Bigr)\gcd(e_1,e_2)\log\lcm(e_1,e_2).
\end{align}
Crucially, these quantities exhibit significant cancellation: all off-diagonal values of $G$ are zero, while the diagonal of $G$ coincides with Euler's $\varphi$ function. Similarly, $L$ vanishes whenever the ratio of its arguments is not a power of a prime, and admits a description in terms of standard arithmetic functions otherwise. The proofs of these assertions are contained in Section \ref{sec:arithm_sums}.

In conclusion, only near-diagonal contributions remain in \eqref{eq:int_var}, and one arrives at the following result. In its statement, a new arithmetic quantity attached to $A$ appears, namely
\begin{equation}
    \label{eq:t_of_a}
     T(A) \vcentcolon= \sum_{d\in D(A)}\frac{\varphi(d)}{m_A(d)}\Bigl(\log\Bigl(\frac{m_A(d)}{d}\Bigr)-S(d)\Bigr),
\end{equation}
where we define
\begin{equation}
    \label{eq:s_of_d}
    S(d) \vcentcolon= \sum_{p \mid d}\frac{\log p}{p-1}
\end{equation}
Here and throughout, unless specified otherwise, `$p$' denotes a prime number.

\begin{proposition}
    \label{prop:int_var}
    Let $A \subseteq \N$ be a finite non-empty set. Then for all sufficiently small $\varepsilon > 0$, we have
    \[
        \int_{\varepsilon}^{\infty}\frac{\Var(W_{A,t})}{\E[{W_{A,t}}]^2}\,dt = \frac{\log(1/\varepsilon)+C_0}{H(A)} + \frac{T(A)}{H(A)^2} + \varepsilon,
    \]
    where $C_0 \vcentcolon= 1 - \log(2\pi)$ is an absolute constant.
\end{proposition}

From here, it is a short step to an upper bound in terms of the quantity $H(A)$.

\begin{theorem}
    \label{thm:bound_in_terms_of_h}
    Let $A \subseteq \N$ be a finite non-empty set. Then
    \[
        \E\bigl[f_A(U)\bigr] \leq \frac{\log H(A) + 2-\log\pi}{H(A)}.
    \]
\end{theorem}

The final ingredient is the following general estimate for $H(A)$, which together with Theorem \ref{thm:bound_in_terms_of_h} yields Theorem \ref{thm:min_fract_part_new}. For a proof, see Section \ref{sec:nt_estimates}.

\begin{proposition}
    \label{prop:bound_on_h}
    For any set $A \subseteq \N$ of size $n$, we have $H(A) \gg \frac{n}{\log n}$.
\end{proposition}

\section{The integrated covariances}\label{sec:int_cov}

The goal of this section is to prove the following.

\begin{proposition}
    \label{prop:int_cov}
    Let $A \subseteq \N$ be a finite non-empty set. Then for all sufficiently small $\varepsilon > 0$, the expression \eqref{eq:int_var} equals exactly
    \[
        \frac{1}{H(A)^2}\sum_{d_1,d_2\in D(A)}\Biggl(\Bigl(\frac{\log(1/\varepsilon)+C_0}{M_A(d_1,d_2)}-\frac{1}{2}J\Bigl(\frac{1}{m_A(d_1)},\frac{1}{m_A(d_2)}\Bigr)\Bigr)G(d_1,d_2) - \frac{L(d_1,d_2)}{M_A(d_1,d_2)}\Biggr) + \varepsilon,
    \]
    where the functions $G$ and $L$ are as defined in \eqref{eq:arithm_sums_g}, \eqref{eq:arithm_sums_l} respectively. Here, for $d_1,d_2 \in D(A)$ we write $M_A(d_1,d_2) \vcentcolon= \max(m_A(d_1),m_A(d_2))$ and for $\alpha,\beta > 0$ we let
    \[
        J(\alpha,\beta) \vcentcolon= \alpha\log\alpha + \beta\log\beta - |\alpha-\beta|\log|\alpha-\beta|.
    \]
\end{proposition}

To prove Proposition \ref{prop:int_cov}, it will be convenient to observe the following identity for $\Delta$:
\begin{equation}
    \label{eq:polarization}
    \Delta(x,y) = \frac{1}{2}\bigl(\psi(x)+\psi(y)-\psi(x-y)\bigr),
\end{equation}
where we define
\[
    \psi \colon \R \to \R, \quad z \mapsto \{z\}(1-\{z\}).
\]
This can be verified by recalling that $\Delta(x,y) = \Cov\bigl(1_{\{U < x\}}, 1_{\{U < y\}}\bigr)$, noting that $\psi$ satisfies $\psi(b-a) = \Var\bigl(1_{\{a \leq U < b\}}\bigr)$ for $0 \leq a \leq b < 1$, and using bilinearity of covariance to expand
\[
    \Var\bigl(1_{\{\min(x,y) \leq U < \max(x,y)\}}\bigr) = \Var\bigl(1_{\{U < \max(x,y)\}} - 1_{\{U < \min(x,y)\}}\bigr).
\]
The identity \eqref{eq:polarization} naturally leads to the following lemma.

\begin{lemma}
    \label{lm:integral}
    Let $a > 0$ be arbitrary and let $\varepsilon \in (0,1/a]$. Then we have
    \[
        \int_{\varepsilon}^{\infty}\frac{\psi(at)}{t^2}\,dt = a\Bigl(\log\Bigl(\frac{1}{\varepsilon}\Bigr) - \log a+\varepsilon a + C_0\Bigr),
    \]
    where $C_0 \vcentcolon= 1 - \log(2\pi)$.
\end{lemma}
\begin{proof}
    By changing variables, we may assume that $a = 1$ and $\varepsilon \in (0,1]$. Hence, since the integrand is non-negative, it suffices to show that
    \[
        \lim_{N\to\infty}\int_{\varepsilon}^{N}\frac{\psi(x)}{x^2}\,dx = \log\Bigl(\frac{1}{\varepsilon}\Bigr) + \varepsilon + C_0.
    \]
    But note that for $N \in \N$ we have
    \begin{align*}
        \int_{\varepsilon}^{N}\frac{\psi(x)}{x^2}\,dx &= \int_{\varepsilon}^{1}\frac{1-x}{x}\,dx + \sum_{k=1}^{N-1}\int_{k}^{k+1}\frac{(x-k)(k+1-x)}{x^2}\,dx\\
        &= \log(1/\varepsilon) + \varepsilon - 1 + \sum_{k=1}^{N-1}\Bigl((2k+1)\log\Bigl(\frac{k+1}{k}\Bigr)-2\Bigr).
    \end{align*}
    Using summation by parts, the last sum can be rewritten as
    \[
        (2N-1)\log N - 2\sum_{k=2}^{N-1}\log k - 2(N - 1) = (2N+1)\log N - 2\log(N!) - 2N + 2.
    \]
    The conclusion now follows from Stirling's approximation $N! \sim \sqrt{2\pi N}(N/e)^N$.
\end{proof}

We are now ready to prove Proposition \ref{prop:int_cov}.

\begin{proof}[Proof of Proposition \ref{prop:int_cov}]
    Using the identity \eqref{eq:polarization} and Lemma \ref{lm:integral}, the integral appearing in \eqref{eq:int_var} can be expressed for $\varepsilon \leq 1/\max A$ as
    \[
        \frac{\lcm(e_1,e_2)}{M_A(d_1,d_2)}\Bigl(\log\Bigl(\frac{1}{\varepsilon}\Bigr)+C_0\Bigr) - \frac{1}{2}J\Bigl(\frac{\lcm(e_1,e_2)}{m_A(d_1)}, \frac{\lcm(e_1,e_2)}{m_A(d_2)}\Bigr) + \frac{\lcm(e_1,e_2)^2}{m_A(d_1)m_A(d_2)}\varepsilon.
    \]
    The latter term accounts for a total contribution of
    \[
        \frac{\varepsilon}{H(A)^2}\sum_{d_1,d_2\in D(A)}\sum_{\substack{e_1\mid d_1\\e_2\mid d_2}}\mu\Bigl(\frac{d_1}{e_1}\Bigr)\mu\Bigl(\frac{d_2}{e_2}\Bigr)\frac{\gcd(e_1,e_2)^2}{e_1e_2}\cdot\frac{\lcm(e_1,e_2)^2}{m_A(d_1)m_A(d_2)}.
    \]
    Note that this equals
    \[
        \frac{\varepsilon}{H(A)^2}\sum_{d_1,d_2\in D(A)}\frac{\Bigl(\sum_{e_1\mid d_1}\mu\Bigl(\frac{d_1}{e_1}\Bigr)e_1\Bigr)\Bigl(\sum_{e_2\mid d_2}\mu\Bigl(\frac{d_2}{e_2}\Bigr)e_2\Bigr)}{m_A(d_1)m_A(d_2)} = \frac{\varepsilon}{H(A)^2}\sum_{d_1,d_2\in D(A)}\frac{\varphi(d_1)\varphi(d_2)}{m_A(d_1)m_A(d_2)} = \varepsilon.
    \]
    The former two terms can be checked to give the correct total contribution by using the identity
    \[
        J\bigl(\kappa\alpha, \kappa\beta\bigr) = \kappa\bigl(J(\alpha,\beta) + 2\min(\alpha,\beta)\log\kappa\bigr)
    \]
    with $\alpha = 1/m_A(d_1)$, $\beta = 1/m_A(d_2)$ and $\kappa = \lcm(e_1,e_2)$. This concludes the proof.
\end{proof}

\section{The arithmetic sums}\label{sec:arithm_sums}

The purpose of this section is to evaluate the arithmetic sums $G$ and $L$ defined in Section \ref{sec:outline}.

\begin{lemma}
    \label{lm:arithm_sums}
    For any $d_1,d_2 \in \N$, we have
    \begin{align*}
        G(d_1,d_2) &= \begin{cases}\varphi(d_1) & \text{if } d_1 = d_2\\0 & \text{else}\end{cases},\\
        L(d_1,d_2) &= \begin{cases}\varphi(d_1)(\log d_1 - S(d_1)) & \text{if } d_1 = d_2\\\varphi(\min(d_1,d_2))\log p & \text{if }\frac{\max(d_1,d_2)}{\min(d_1,d_2)} = p^j \text{ for some } p \text{ and } j \in \N\\0 & \text{else}\end{cases}.
    \end{align*}
\end{lemma}
\begin{proof}
    Given $d_1,d_2 \in \N$, consider the following real Dirichlet polynomial:
    \[
        K_{d_1,d_2}(x) = \sum_{\substack{e_1\mid d_1\\e_2\mid d_2}}\mu\Bigl(\frac{d_1}{e_1}\Bigr)\mu\Bigl(\frac{d_2}{e_2}\Bigr)\gcd(e_1,e_2)\lcm(e_1,e_2)^x.
    \]
    We may interpret the quantities of interest in terms of $K_{d_1,d_2}$ as $G(d_1,d_2) = K_{d_1,d_2}(0)$ and $L(d_1,d_2) = K_{d_1,d_2}'(0)$. To evaluate these, note that for any fixed $x$, the function $(d_1,d_2) \mapsto K_{d_1,d_2}(x)$ can be viewed as the Dirichlet convolution of the two-variable functions
    \[
        (m,n) \mapsto \mu(m)\mu(n), \quad (m,n) \mapsto \gcd(m,n)\lcm(m,n)^x,
    \]
    both of which are multiplicative. As such, this function is itself multiplicative, and hence admits the product representation
    \[
        K_{d_1,d_2}(x) = \prod_{p}K_{p^{\nu_p(d_1)},p^{\nu_p(d_2)}}(x),
    \]
    where $\nu_p$ denotes $p$-adic valuation (for the precise definition of multivariate multiplicative functions and their relevant properties, see e.g.\ \cite[\S3.1,\S4.1]{toth}). Given any prime $p$, it is routine to verify that for any $\alpha_1,\alpha_2 \in \N_0$ we have
    \[
        K_{p^{\alpha_1},p^{\alpha_2}}(x) = \sum_{\substack{0 \leq j_1 \leq \alpha_1\\0\leq j_2 \leq \alpha_2}}\mu(p^{\alpha_1-j_1})\mu(p^{\alpha_2-j_2})p^{\min(j_1,j_2)+x\max(j_1,j_2)}
    \]
    and hence that
    \[
        K_{p^{\alpha_1},p^{\alpha_2}}(x) = \begin{cases}\varphi(p^{\min(\alpha_1,\alpha_2)})p^{x(\max(\alpha_1,\alpha_2)-1)}(p^x-1) & \text{if } \alpha_1 \neq \alpha_2\\p^{(\alpha_1-1)(x+1)}(p^{x+1}-2p^x+1) & \text{if } \alpha_1 = \alpha_2 > 0\\1 & \text{if } \alpha_1 = \alpha_2 = 0\end{cases}.
    \]
    In particular, it is now clear that
    \[
        K_{p^{\alpha_1},p^{\alpha_2}}(0) = \begin{cases}0 & \text{if } \alpha_1 \neq \alpha_2\\\varphi(p^{\alpha_1}) & \text{if } \alpha_1 = \alpha_2\end{cases},
    \]
    and hence the claim about the value of $G(d_1,d_2)$ follows from the multiplicativity of $\varphi$. To verify the claim about $L(d_1,d_2)$, note that the product rule gives
    \[
        K_{d_1,d_2}'(0) = \sum_{p}K_{p^{\nu_p(d_1)},p^{\nu_p(d_2)}}'(0)\prod_{q \neq p}K_{q^{\nu_q(d_1)},q^{\nu_q(d_2)}}(0),
    \]
    where the product is over all primes $q \neq p$. In the above sum, the only surviving terms are those for which $\nu_q(d_1) = \nu_q(d_2)$ for all $q \neq p$. In particular, if the exponents in the prime factorizations of $d_1$ and $d_2$ differ for at least two primes, the sum is zero. Otherwise, if they are different for exactly one prime $p$, we obtain
    \[
        \varphi\bigl(p^{\min(\nu_p(d_1),\nu_p(d_2))}\bigr)\log p\cdot\prod_{q\neq p}\varphi(q^{\nu_q(d_1)}) = \varphi\bigl(\min(d_1,d_2)\bigr)\log p,
    \]
    as required. Finally, if $d_1 = d_2$, then it's perhaps simplest to consider logarithmic derivatives:
    \[
        \bigl(\log K_{d_1,d_2}\bigr)'(0) = \sum_{p}\bigl(\log K_{p^{\nu_p(d_1)},p^{\nu_p(d_2)}}\bigr)'(0) = \sum_{p\mid d_1}\Bigl(\nu_p(d_1)\log p - \frac{\log p}{p-1}\Bigr) = \log d_1 - S(d_1),
    \]
    whence multiplying by $K_{d_1,d_2}(0) = \varphi(d_1)$ gives the claim.
\end{proof}

\section{Number-theoretic estimates}\label{sec:nt_estimates}

Our first result in this section is the following refinement of Proposition \ref{prop:bound_on_h}.

\begin{proposition}
    \label{prop:bound_on_h_explicit}
    For $j \geq 1$, let $p_j$ denote the $j$-th prime. For $n \geq 0$, define
    \[
        P_n \vcentcolon= \prod_{j=1}^{n}\Bigl(1-\frac{1}{p_j}\Bigr),
    \]
    and for $n \geq 1$, write $R_n \vcentcolon= \sum_{j=0}^{n-1}P_j$. Then for any $A \subseteq \N$ of size $n$, we have
    \[
        H(A) \geq R_n.
    \]
    Moreover, $R_n \sim \frac{e^{-\gamma}n}{\log n}$, where $\gamma$ is the Euler--Mascheroni constant, so in particular $H(A) \gg \frac{n}{\log n}$.
\end{proposition}

We also give an absolute upper bound on the quantity $T(A)$ introduced in Section \ref{sec:outline}.

\begin{lemma}
    \label{lm:t_is_non_positive}
    For any finite $A \subseteq \N$, we have $T(A) \leq 0$.
\end{lemma}

In the proofs of the above results, the following device will prove useful. For a positive integer $a$, let $X_a$ be a random variable taking values in the set of divisors of $a$ according to the distribution
\[
    \P(X_a = d) = \frac{\varphi(d)}{a}.
\]
Observe that the exponents in the prime factorization of $X_a$ are independent and follow the distribution
\[
    \P(\nu_p(X_a) = j) = \frac{\varphi(p^j)}{p^{\nu_p(a)}} 
\]
for $p \mid a$ and $0 \leq j \leq \nu_p(a)$.

Given a finite set $A \subseteq \N$, enumerate it as $A = \{a_1,\ldots,a_n\}$ where $a_1 < \ldots < a_n$. We partition $D(A)$ according to the value of $m_A$ by defining for $i \in [n]$ the set
\[
    B_i \vcentcolon= \{d \in D(A) : m_A(d) = a_i\} = \bigl\{d \in \N : d \mid a_i,\ \forall j < i\ d\nmid a_j\bigr\}.
\]
We are now ready to proceed with the main proofs in this section.

\begin{proof}[Proof of Proposition \ref{prop:bound_on_h_explicit}]
    By the above discussion, we have
    \begin{equation}
        \label{eq:h_exp}
        H(A) = \sum_{i=1}^{n}\sum_{d\in B_i}\frac{\varphi(d)}{a_i} = \sum_{i=1}^{n}\P(X_{a_i} \in B_i) = \sum_{i=1}^{n}\P\bigl(\forall j < i\ X_{a_i} \nmid a_j\bigr).
    \end{equation}
    Choosing for each $j < i$ a prime $q_{i,j}$ such that $\nu_{q_{i,j}}(a_i) > \nu_{q_{i,j}}(a_j)$ and setting 
    \[
        Q_i \vcentcolon= \bigl\{q_{i,j} : j \in \{1,\ldots,i-1\}\bigr\},
    \]
    the above probability can be lower-bounded by
    \begin{align*}
        \P\bigl(\forall j \in \{1,\ldots,i-1\}\ \nu_{q_{i,j}}(X_{a_i}) = \nu_{q_{i,j}}(a_i)\bigr) = \prod_{q \in Q_i}\P\bigl(\nu_q(X_{a_i}) = \nu_q(a_i)\bigr) = \prod_{q \in Q_i}\Bigl(1-\frac{1}{q}\Bigr) \geq P_{i-1}.
    \end{align*}
    The bound $H(A) \geq R_n$ now follows by summing over all $i \in [n]$.

    To conclude the proof, it remains to note that by Mertens' third theorem \cite[Theorem I.1.11]{tenenbaum} (and Chebyshev's estimate \cite[Theorem I.1.3]{tenenbaum}, by which $p_j \asymp j\log j$), we have $P_j \sim e^{-\gamma}/\log j$. The desired asymptotic for $R_n$ now follows by summing over all $j \in \{0,\ldots,n-1\}$.
\end{proof}

\begin{proof}[Proof of Lemma \ref{lm:t_is_non_positive}]
    Note that we have
    \[
        T(A) = \sum_{i=1}^{n}\sum_{d \in B_i}\frac{\varphi(d)}{a_i}\Bigl(\log\Bigl(\frac{a_i}{d}\Bigr)-S(d)\Bigr) = \sum_{i=1}^{n}\E\bigl[V_{a_i}(X_{a_i})1_{B_i}(X_{a_i})\bigr],
    \]
    where for $a \in \N$ we define the function
    \[
        V_a \colon \N \to \R,\quad d \mapsto \log(a/d)-S(d).
    \]
    The key point is that both $-V_{a_i}(X_{a_i})$ and $1_{B_i}(X_{a_i})$ are increasing functions of the exponents in the factorization of $X_{a_i}$ into primes. Since these exponents are independent random variables, \cite[Theorem 2.1]{esary-proschan-walkup} implies that they are associated, so in particular
    \[
        \Cov\bigl(-V_{a_i}(X_{a_i}), 1_{B_i}(X_{a_i})\Bigr) \geq 0.
    \]
    But $V_{a_i}(X_{a_i})$ has mean zero: both $\E\bigl[\log(a_i/X_{a_i})\bigr]$ and $\E\bigl[S(X_{a_i})\bigr]$ are equal to
    \[
        \sum_{p}\bigl(1-p^{-\nu_p(a_i)}\bigr)\frac{\log p}{p-1}.
    \]
    We thus obtain that
    \[
        \E\bigl[V_{a_i}(X_{a_i})1_{B_i}(X_{a_i})\bigr] \leq 0,
    \]
    from which the conclusion follows by summing over all $i \in [n]$.
\end{proof}

\section{Putting everything together}\label{sec:pf_of_thm}

In this section we put together everything we have developed so far and tie up some loose ends from Section \ref{sec:outline}. We first explain how to modify the definition of the random variable $W_{A,t}$ when $t \geq 1$ so that the identity \eqref{eq:var_identity} remains valid. The idea is to let
\[
    W_{A,t} = \sum_{d \in D(A)}\chi_{d,t/m_A(d)}^*(U)
\]
with $\chi_{d,\delta}^*$ defined as in \eqref{eq:chi}. It is clear that, for $t \in [0,1)$, this agrees with the previous definition \eqref{eq:counting_coprime}. The expansion \eqref{eq:var_expansion} goes through; however, in order to calculate
\[
    \Cov\bigl(\chi_{e_1,t/m_A(d_1)}(U), \chi_{e_2,t/m_A(d_2)}(U)\bigr),
\]
we may no longer apply Lemma \ref{lm:covariance} directly. Nevertheless, there is a simple workaround for this issue, which is to note that
\[
    \chi_{a,\alpha}(U) = \lfloor a\alpha\rfloor + 1_{\{\{aU\} < \{a\alpha\}\}} = \lfloor a\alpha\rfloor + 1_{E_{a,\{a\alpha\}}}(U).
\]
The constant terms do not affect the covariance, and an application of Lemma \ref{lm:covariance} recovers \eqref{eq:var_identity}.

Even though we were able to use Lemma \ref{lm:covariance} as a black box, we give a separate proof of the exact statement we need in Appendix \ref{app:covariance}. Besides making the paper self-contained, this alternative perspective may be enlightening in its own right.

We are now in a position to prove Proposition \ref{prop:int_var}.

\begin{proof}[Proof of Proposition \ref{prop:int_var}]
    Using Proposition \ref{prop:int_cov} and the first part of Lemma \ref{lm:arithm_sums}, the left-hand side can be written as
    \begin{equation}
        \label{eq:simplified}
        \frac{1}{H(A)^2}\Biggl(\sum_{d\in D(A)}\frac{\varphi(d)}{m_A(d)}\bigl(\log(1/\varepsilon)+C_0+\log m_A(d)\bigr) - \sum_{d_1,d_2\in D(A)}\frac{L(d_1,d_2)}{M_A(d_1,d_2)}\Biggr) + \varepsilon.
    \end{equation}
    Using the second part of Lemma \ref{lm:arithm_sums} and noting that for any $d \in D(A)$, all its divisors are also present in $D(A)$, the second sum can be expressed as
    \[
        \sum_{d\in D(A)}\frac{\varphi(d)}{m_A(d)}(\log d - S(d)) + 2\sum_{d\in D(A)}\frac{1}{m_A(d)}\sum_{p\mid d}\sum_{j=1}^{\nu_p(d)}\varphi\bigl(d/p^j\bigr)\log p.
    \]
    The second term in the above expression simplifies to
    \[
        2\sum_{d\in D(A)}\frac{\varphi(d)}{m_A(d)}S(d).
    \]
    Hence, \eqref{eq:simplified} becomes
    \[
        \frac{1}{H(A)^2}\sum_{d\in D(A)}\frac{\varphi(d)}{m_A(d)}\bigl(\log(1/\varepsilon) + C_0 + \log m_A(d) - \log d - S(d)\bigr) + \varepsilon,
    \]
    whereupon the conclusion follows by pulling out the constant terms and recalling \eqref{eq:t_of_a}.
\end{proof}

We can now prove Theorem \ref{thm:bound_in_terms_of_h}.

\begin{proof}[Proof of Theorem \ref{thm:bound_in_terms_of_h}]
    Recall that we have
    \begin{equation}
        \label{eq:exp_bound}
        \E\bigl[f_A(U)\bigr] = \int_{0}^{1}\P(W_{A,t} = 0)\,dt \leq \eta + \int_{\eta}^{\infty}\frac{\Var(W_{A,t})}{\E[W_{A,t}]^2}\,dt,
    \end{equation}
    where $\eta \in (0,1]$ is any parameter of our choice. We would ideally want to apply this with $\eta = \varepsilon$ for the $\varepsilon$ occurring in Proposition \ref{prop:int_var}. However, Proposition \ref{prop:int_var} applies only for $\varepsilon$ small enough, which might be just too small for the ensuing bound to be useful. Nevertheless, we can remove the `sufficiently small' qualification on $\varepsilon$ using the following trick. Given $\eta$ as above, pick $\varepsilon \in (0,\eta]$ sufficiently small so that Proposition \ref{prop:int_var} applies. By combining its conclusion with Lemma \ref{lm:t_is_non_positive} to eliminate the second term on the right-hand side, we obtain
    \[
        \int_{\varepsilon}^{\infty}\frac{\Var(W_{A,t})}{\E[W_{A,t}]^2}\,dt \leq \frac{\log(1/\varepsilon) + C_0}{H(A)} + \varepsilon.
    \]
    On the other hand, since $W_{A,t}$ takes values in the non-negative integers, we have ${W_{A,t}}^2 \geq W_{A,t}$ deterministically and hence
    \[
        \int_{\varepsilon}^{\eta}\frac{\Var(W_{A,t})}{\E[W_{A,t}]^2}\,dt \geq \int_{\varepsilon}^{\eta}\frac{\E[W_{A,t}]-\E[W_{A,t}]^2}{\E[W_{A,t}]^2}\,dt = \int_{\varepsilon}^{\eta}\Bigl(\frac{1}{\E[W_{A,t}]}-1\Bigr)\,dt = \frac{\log(\eta/\varepsilon)}{H(A)}-\eta+\varepsilon.
    \]
    Subtracting this from the preceding display then gives
    \[
        \int_{\eta}^{\infty}\frac{\Var(W_{A,t})}{\E[W_{A,t}]^2}\,dt \leq \frac{\log(1/\eta)+C_0}{H(A)} + \eta,
    \]
    as desired. Combining this with \eqref{eq:exp_bound} now gives
    \[
        \E\bigl[f_A(U)\bigr] \leq 2\eta + \frac{\log(1/\eta)+C_0}{H(A)}.
    \]
    It remains only to optimize this bound in $\eta$. Choosing $\eta = 1/(2H(A))$, which is permissible by Proposition \ref{prop:bound_on_h_explicit}, we obtain the required bound.
\end{proof}

At last, we can give the short deduction of Theorem \ref{thm:min_fract_part_new}.

\begin{proof}[Proof of Theorem \ref{thm:min_fract_part_new}]
    As explained in the proof outline of Theorem \ref{thm:min_fract_part_old}, we may assume without loss of generality that $A \subseteq \N$. By Theorem \ref{thm:bound_in_terms_of_h}, we know that $\E\bigl[f_A(U)\bigr] \leq \vartheta(H(A))$, where
    \[
        \vartheta \colon (0,\infty) \to \R, \quad h \mapsto \frac{\log h + 2 - \log\pi}{h}.
    \]
    We have $\vartheta'(h) = (\log\pi-1-\log h)/h^2$, so $\vartheta$ is decreasing on $[\pi/e,\infty)$. For $n \geq 2$, Proposition \ref{prop:bound_on_h_explicit} thus implies
    \[
        \vartheta(H(A)) \leq \vartheta(R_n) \ll \frac{(\log n)^2}{n},
    \]
    as required.
\end{proof}

\begin{remark}
    \label{rem:constants}
    When all elements of $A$ have the same sign, the proof of Theorem \ref{thm:min_fract_part_new} delivers $e^{\gamma} + o(1) \approx 1.7811$ as the value of the implied constant in \eqref{eq:min_frac_part_new}. When we allow both positive and negative elements, we get a constant which is at most twice as large. In fact, it is possible to carry out the proof so as to recover the same constant in the general case. However, this entails extra difficulties regarding notation, and the resulting bound is probably not sharp anyway. For these reasons, we have opted for the current presentation.
\end{remark}

\section{Concluding remarks}\label{sec:conc_rem}

In this paper, we made an improvement to the upper bound on $g(m)$, the minimum common size of $A_1,\ldots,A_m \subseteq \Z$ guaranteeing the existence of a non-vanishing transversal. We achieved this by establishing an improved upper bound on the expected value of $f_A(U)$, the minimum fractional part in a uniformly random dilate of a set $A \subseteq \Z\setminus\{0\}$. The crucial idea behind this result was to find a more efficient representation of the sublevel sets of the random variable $f_A(U)$ as a union of intervals. This reduced the overlap between the intervals constituting this union, which in general resulted in a smaller variance in the count of intervals containing $U$. At the same time, it decreased the expectation of the same count, making it proportional to the quantity $H(A)$ instead of $|A|$. This approach naturally led to a bound on the expectation of $f_A(U)$ in terms of $H(A)$, namely Theorem \ref{thm:bound_in_terms_of_h}.

Even though $H(A)$ in general does not exceed $|A|$ (as can be seen for example from \eqref{eq:h_exp}), we've shown in Proposition \ref{prop:bound_on_h_explicit} that it is always smaller by at most a logarithmic factor. It is natural to wonder whether this bound is sharp, since a sharpening would immediately lead to an improvement in Theorem \ref{thm:min_fract_part_new}. The following result confirms that, in fact, the bound in Proposition \ref{prop:bound_on_h_explicit} is best possible up to lower-order terms. Its proof is given in Appendix \ref{app:small_h}.

\begin{proposition}
    \label{prop:small_h}
    For any $n \in \N$, there exists a set $A \subseteq \N$ of size $n$ such that 
    \[
        H(A) \leq \frac{n}{(\log n)^{1-o(1)}}.
    \]
\end{proposition}

Another natural question is whether the bound in Theorem \ref{thm:min_fract_part_new} is sharp. The following result shows that it is at most one logarithmic factor away from the truth.

\begin{theorem}
    \label{thm:initial_segment}
    As $n \to \infty$, we have the asymptotic
    \[
        \E\bigl[f_{[n]}(U)\bigr] \sim \frac{3\log n}{\pi^2n}.
    \]
\end{theorem}

Theorem \ref{thm:initial_segment} is a consequence of a result of Broder, Charikar, Frieze and Mitzenmacher on min-wise independent permutations \cite[Theorem 11]{broder-charikar-frieze-mitzenmacher}. In fact, the desired estimate can be extracted from the proof of this result: the quantity $A_1$ defined in the course of the proof corresponds precisely to the above expectation. The proof supplies the desired asymptotic, with a multiplicative error of the form $1 + O(1/\log n)$. 

Theorem \ref{thm:initial_segment} also rules out the possibility of proving the Goddyn--Tarsi conjecture by the approach taken in the present paper. Nevertheless, it might still be possible to improve Theorem \ref{thm:min_fract_part_new}; we are not aware of any choice of $A$ that produces a significantly larger expected value of $f_A(U)$ than $A = [n]$. Hence, we ask the following question.

\begin{question}
    \label{q:max_exp_fract_part}
    For a given positive integer $n$, what is the supremum of $\E\bigl[f_A(U)\bigr]$ as $A$ ranges over all $n$-element subsets of $\N$?
\end{question}


Finally, it would be remiss not to mention the connection of the problems from this paper to the lonely runner conjecture. This conjecture is originally due to Wills \cite{wills}, and independently Cusick \cite{cusick}; its present name was incidentally given by Goddyn \cite{bienia-goddyn-gvozdjak-sebo-tarsi}. It states that, for any set $A$ of $n$ positive integers,
\[
    \max_{\theta \in \T}\min_{a\in A}\,\lVert a\theta\rVert_{\T} \geq \frac{1}{n+1},
\]
where $\lVert\cdot\rVert_{\T}$ denotes the distance to the nearest integer. Changing the maximum to an average and the toroidal distance to a fractional part, we obtain exactly the quantity from Theorem \ref{thm:main_ingredient}.

The trivial lower bound of $\frac{1}{2n}$ for the lonely runner conjecture follows by observing that for any $a \in A$ the set of all $\theta \in \T$ such that $\lVert a\theta\rVert_{\T} < \delta$ has measure $2\delta$ and hence by the union bound, the set of all $\theta$ where this holds for at least one $a \in A$ has measure at most $2\delta n$. In particular, if $\delta < \frac{1}{2n}$, then there exists $\theta \in \T$ such that $\lVert a\theta\rVert_{\T} \geq \delta$ for all $a \in A$.

Essentially the same argument as above gives a matching lower bound on $\E[f_A(U)]$ for any $n$-element set $A$. Indeed, the set of all $\theta \in [0,1)$ for which there exists $a \in A$ with $\{a\theta\} < \frac{1}{2n}$ has measure at most $\frac{1}{2}$. Thus, for all $\theta$ in a set of measure at least $\frac{1}{2}$, we have $\{a\theta\} \geq \frac{1}{2n}$ for all $a \in A$, and hence $\E[f_A(U)] \geq \frac{1}{4n}$. This bound is tight up to constants: for a generic\footnote{Meaning the elements are selected uniformly at random from a long initial segment of $\N$.} choice of $A$, the expected value of $\E[f_A(U)]$ is asymptotically $\frac{1}{n}$.

Recently there have been important advances on the lonely runner conjecture. In 2025, Bedert \cite{bedert} proved a lower bound of $\frac{1}{2n} + n^{-5/3-o(1)}$, improving on previous bounds of the form $\frac{1}{2n} + n^{-2+o(1)}$. In a very recent preprint \cite{hu-jing}, Hu and Jing obtained a further improvement by establishing a bound of the form $\frac{1+c}{2n}$, where $c > 0$ is an absolute constant. Their arguments seem to have certain similarities with the methods of the present paper, in particular in the way they deal with the overlap between different intervals in unions of the type \eqref{eq:interval_union}. It remains an interesting problem to investigate whether this connection can be further exploited in either direction.

\section{AI statement}\label{sec:ai}

All new mathematical arguments in the paper are due to GPT-6 Astra, prompted by the author. In particular, the key new result, Theorem \ref{thm:main_ingredient}, was established by GPT-6 Astra. The proofs were digested and rewritten by the author, who takes full responsibility for the content of the paper.

On October 2nd, 2026, the author asked Astra to give a best possible upper bound on the probability appearing in \eqref{eq:layer_cake_first}. Astra gave the bound obtained by inputting \eqref{eq:gal_sums} into the Alon--Ruzsa argument as described in Section \ref{sec:outline}. On October 3rd, the author prompted Astra to improve this bound, pointing it to an approach along the lines of the recent breakthrough of Bedert \cite{bedert} on the lonely runner conjecture. In response, Astra was able to obtain an improvement in a regime where $t$ is small in terms of $n$. However, when asked how this bound is reflected in the expected value of $f_A(U)$, it replied that it yields no improvement whatsoever. The author then asked it to improve the bound on the expectation of $f_A(U)$ itself, which resulted in Theorem \ref{thm:main_ingredient}. Subsequently, while writing this paper, the author repeatedly asked Astra to explain its approach and replace certain arguments by explicit references to the literature. Finally, GPT-6.2 Astra was used to proofread the paper.

\appendix

\section{Covariance of dilates modulo \texorpdfstring{$1$}{1}}\label{app:covariance}

In this section we give a self-contained proof of the following slight generalization of Lemma \ref{lm:covariance}. Our argument is different from the original one of Alon and Ruzsa \cite{alon-ruzsa}.

\begin{lemma}
    \label{lm:covariance_general}
    Let $a_1,a_2 \in \N$ and $\alpha_1,\alpha_2 \geq 0$. Write $\ell \vcentcolon= \lcm(a_1,a_2)$ and $a_i' \vcentcolon= a_i/\gcd(a_1,a_2)$ for $i \in [2]$. Then we have
    \[
        \Cov\bigl(\chi_{a_1,\alpha_1}(U), \chi_{a_2,\alpha_2}(U)\bigr) = \frac{\Delta\bigl(\{\ell\alpha_1\}, \{\ell\alpha_2\}\bigr)}{a_1'a_2'}.
    \]
\end{lemma}
\begin{proof}
    Write $U = (J + U')/\ell$, where $J$ and $U'$ are independent random variables uniformly distributed in $\{0,1,\ldots,\ell-1\}$ and $[0,1)$ respectively. We first note that, for any fixed $i \in [2]$ and $u \in [0,1)$, the function
    \[
        \Z \to \Z, \quad j \mapsto \chi_{a_i,\alpha_i}\bigl((j+u)/\ell\bigr)
    \]
    is $\ell/a_i$-periodic. Thus, $\chi_{a_i,\alpha_i}((J+u)/\ell)$ is a function of $J+(\ell/a_i)\Z$. But the latter is uniformly distributed in $\Z/(\ell/a_i)\Z$, so we have
    \[
        \E\Bigl[\chi_{a_i,\alpha_i}\bigl((J+u)/\ell\bigr)\Bigr] = \frac{1}{\ell/a_i}\sum_{0 \leq r < \ell/a_i}\chi_{a_i,\alpha_i}\bigl((r+u)/\ell\bigr).
    \]
    Expanding out the definition of $\chi_{a_i,\alpha_i}$, we obtain that this equals
    \[
        \frac{1}{\ell/a_i}\sum_{\substack{0\leq r < \ell/a_i\\j \in \Z}}1_{[j/a_i,j/a_i+\alpha_i)}\bigl((r+u)/\ell\bigr)
        = \frac{1}{\ell/a_i}\sum_{\substack{0\leq r < \ell/a_i\\j \in \Z}}1_{[j\ell/a_i-r, j\ell/a_i-r + \ell\alpha_i)}(u)
        = \frac{1}{\ell/a_i}\chi_{1,\ell\alpha_i}(u),
    \]
    where we used in the last equality that the map
    \[
        \bigl\{0,1,\ldots,\ell/a_i-1\bigr\}\times\Z \to \Z, \quad (r,j) \mapsto j\ell/a_i-r
    \]
    is a bijection.

    Furthermore, by the Chinese remainder theorem, $J+(\ell/a_1)\Z$ and $J+(\ell/a_2)\Z$ are independent, so it follows that
    \[
        \E\Biggl[\prod_{i=1}^{2}\chi_{a_i,\alpha_i}\bigl((J+u)/\ell\bigr)\Biggr] = \prod_{i=1}^{2}\E\Bigl[\chi_{a_i,\alpha_i}\bigl((J+u)/\ell\bigr)\Bigr] = \frac{\chi_{1,\ell\alpha_1}(u)\chi_{1,\ell\alpha_2}(u)}{a_1'a_2'}.
    \]
    Since $J$ and $U'$ are independent, integrating over $u$ yields
    \[
        \E\bigl[\chi_{a_1,\alpha_1}(U)\chi_{a_2,\alpha_2}(U)\bigr] =\frac{1}{a_1'a_2'}\int_{0}^{1}\chi_{1,\ell\alpha_1}(u)\chi_{1,\ell\alpha_2}(u)\,du
        = \frac{\E\bigl[\chi_{1,\ell\alpha_1}(U')\chi_{1,\ell\alpha_2}(U')\bigr]}{a_1'a_2'}.
    \]
    Similarly, for $i \in [2]$ we have
    \[
        \E\bigl[\chi_{a_i,\alpha_i}(U)\bigr] = \int_{0}^{1}\E\Bigl[\chi_{a_i,\alpha_i}\bigl((J+u)/\ell\bigr)\Bigr]\,du = \frac{1}{\ell/a_i}\int_{0}^{1}\chi_{1,\ell\alpha_i}(u)\,du = \frac{\E\bigl[\chi_{1,\ell\alpha_i}(U')\bigr]}{\ell/a_i},
    \]
    so we obtain that
    \[
         \Cov\bigl(\chi_{a_1,\alpha_1}(U), \chi_{a_2,\alpha_2}(U)\bigr) = \frac{\Cov\bigl(\chi_{1,\ell\alpha_1}(U'), \chi_{1,\ell\alpha_2}(U')\bigr)}{a_1'a_2'}.
    \]
    The desired conclusion now follows upon observing that $\chi_{1,\alpha}(U') = \lfloor\alpha\rfloor + 1_{\{U'<\{\alpha\}\}}$ and recalling that $\Cov(1_{\{U'<x\}},1_{\{U'<y\}}) = \Delta(x,y)$ for $x,y \in [0,1)$.
\end{proof}

\section{Sets with small \texorpdfstring{$H(A)$}{H(A)}}\label{app:small_h}

The purpose of this section is to prove Proposition \ref{prop:small_h}, which provides sets $A$ with small $H(A)$ relative to their size. The starting point for the construction is the following lemma.

\begin{lemma}
    \label{lm:inverses}
    Let $B \subseteq \N$ be a finite non-empty set and let $Q$ be a common multiple of its elements. Then for $A = \{Q / b : b \in B\}$, we have
    \begin{equation}
        \label{eq:exp_for_h}
        H(A) = \E\bigl[\max(\{b \in B : b \mid J\} \cup \{0\})\bigr],
    \end{equation}
    where $J$ is a random variable with uniform distribution in $[Q]$. In particular, if $B \subseteq (x,2x]$,
    \begin{equation}
        \label{eq:h_prob}
        H(A) \asymp x\P(\exists b \in B\ b \mid J).
    \end{equation}
\end{lemma}
\begin{proof}
    For any $b \in B$, we have $b \mid Q$, so $b \mid J$ is equivalent to $b \mid \gcd(J, Q)$. The latter is in turn equivalent to $Q/\gcd(J, Q) \mid Q/b$. Thus, the maximum $b \in B$ such that $b \mid J$ equals $Q$ divided by the minimum $a \in A$ such that $Q/\gcd(J,Q) \mid a$. Moreover, the maximum exists if and only if the minimum exists, and the latter happens precisely if $Q/\gcd(J,Q) \in D(A)$. Therefore, the right-hand side of \eqref{eq:exp_for_h} equals
    \[
        \sum_{d \in D(A)}\P\bigl(Q/\gcd(J,Q) = d\bigr)\cdot\frac{Q}{m_A(d)} = \sum_{d\in D(A)}\frac{\varphi(d)}{Q}\cdot\frac{Q}{m_A(d)} = H(A),
    \]
    as required. The second statement now follows since the random variable on the right-hand side of \eqref{eq:exp_for_h} is bounded above and below by constant multiples of $x1_{\{\exists b \in B\ b \mid J\}}$.
\end{proof}

The idea now is to take the set $B$ in Lemma \ref{lm:inverses} to consist of numbers in a dyadic interval $(x,2x]$ having unusually many prime factors. Specifically, we will choose $B$ to be an $n$-element subset of the set
\[
    \mathcal{B}_{x,c} \vcentcolon= \bigl\{b \in \N : x < b \leq 2x,\ \omega(b) \geq c\log\log x\bigr\},
\]
where $c > 1$ is some constant and $\omega(b)$ denotes the number of distinct prime factors of $b$. For such a choice of $B$, the following lemma provides a bound on the probability appearing on the right-hand side of \eqref{eq:h_prob}. In what follows, an important role will be played by the large-deviation rate function
\[
    I \colon (0,\infty) \to \R, \quad u \mapsto u\log u - u + 1.
\]

\begin{lemma}
    \label{lm:div_prob}
    Let $c_1 > 1$ and $c_2 \in (0,1)$ be fixed. Then for $B \subseteq \mathcal{B}_{x,c_1}$, $Q$ a common multiple of the elements of $B$ and $J$ uniformly distributed in $[Q]$, we have
    \[
        \P(\exists b \in B\ b \mid J) \leq (\log x)^{-I(c_1+c_2)+o(1)} + \frac{|B|}{x}(\log x)^{-I(c_2)+o(1)}.
    \]
\end{lemma}
\begin{proof}
    The probability on the left-hand side does not depend on the particular choice of $Q$, so we may assume without loss of generality that $Q = \lcm(1,\ldots,\lfloor2x\rfloor)$. Consider for each prime $p$ the random variable $Y_p \vcentcolon=1_{\{p \mid J\}}$. By the Chinese remainder theorem, $(Y_p)_{p\leq 2x}$ are independent Bernoulli random variables, where $Y_p$ has mean $1/p$. Letting
    \[
        C \vcentcolon= \sum_{p \leq 2x}Y_p,
    \]
    linearity of expectation and Mertens' second theorem \cite[Theorem I.1.9]{tenenbaum} imply
    \[
        \E[C] = \sum_{p \leq 2x}\frac{1}{p} = \log\log x + O(1).
    \]
    Hence, by the Chernoff (upper-tail) bound (see e.g.\ \cite[Theorem 4.4]{mitzenmacher-upfal}), we have
    \[
        \P(C \geq (c_1+c_2)\log\log x) \leq (\log x)^{-I(c_1+c_2)+o(1)}.
    \]
    Returning to our main task, note that the union bound implies
    \[
        \P(\exists b \in B\ b \mid J) \leq \P(C \geq (c_1+c_2)\log\log x) + \sum_{b\in B}\P(b \mid J,\ C < (c_1+c_2)\log\log x).
    \]
    By the above, it suffices to upper bound each summand by $(\log x)^{-I(c_2)+o(1)}/x$. To this end, note that fixing $b \in B$ and writing $\mathcal{E}$ for the event that $b\mid J$, we have $\P(\mathcal{E}) = 1/b < 1/x$. Furthermore, since $\omega(b) \geq c_1\log\log x$,
    \[
        \P\bigl(C < (c_1+c_2)\log\log x\ \big|\ \mathcal{E}\bigr) \leq \P\Bigl(\sum_{p \leq 2x,\ p \nmid b}Y_p < c_2\log\log x\Bigr).
    \]
    But note that
    \[
        \E\Bigl[\sum_{p \leq 2x,\ p \nmid b}Y_p\Bigr] = \sum_{p \leq 2x,\ p \nmid b}\frac{1}{p} = \E[C] - \sum_{p\mid b}\frac{1}{p} \geq \log\log x - \log\log\log x - O(1),
    \]
    where we used the standard fact that
    \[
        \sum_{p\mid b}\frac{1}{p} \leq \log\log\log b + O(1);
    \]
    for a proof, see e.g.\ \cite[Corollary 2.7]{beker-renyi-entropy}. Hence, the desired conclusion follows by applying the Chernoff (lower-tail) bound (see e.g.\ \cite[Theorem 4.5]{mitzenmacher-upfal}).
\end{proof}

We can now give the proof of Proposition \ref{prop:small_h}.

\begin{proof}[Proof of Proposition \ref{prop:small_h}]
    Let $\varepsilon > 0$ be given; we will show that for large enough $n$ there is a set $A \subseteq \N$ of size $n$ such that
    \[
        H(A) \leq \frac{n}{(\log n)^{1-\varepsilon}}.
    \]
    Let $x = n(\log n)^{I(c_1)+1}$, where $c_1 > 1$ is a constant parameter to be determined later. Then by \cite[Theorem II.6.4]{tenenbaum} and Stirling's approximation, for $k = \lceil c_1\log\log x\rceil$ we have
    \[
        |\mathcal{B}_{x,c_1}| \geq \pi_k(2x) - \pi_k(x) = x(\log x)^{-I(c_1)+o(1)} = n(\log n)^{1-o(1)},
    \]
    where $\pi_k(y)$ denotes the number of $b \leq y$ with $\omega(b) = k$. In particular, if $n$ is large enough, we may select $B$ to be a subset of $\mathcal{B}_{x,c_1}$ of size $n$. Letting $Q = \lcm(B)$, $J$ be uniformly distributed in $[Q]$ and $A = \{Q/b : b \in B\}$, Lemmas \ref{lm:inverses} and \ref{lm:div_prob} give
    \begin{align*}
        \frac{H(A)}{n} &\asymp \frac{x}{n}\P(\exists b \in B\ b \mid J)\\
        &\leq \frac{x}{n}(\log x)^{-I(c_1+c_2)+o(1)} + (\log x)^{-I(c_2)+o(1)}\\
        &= (\log n)^{I(c_1)-I(c_1+c_2)+1+o(1)} + (\log n)^{-I(c_2)+o(1)},
    \end{align*}
    where $c_2 \in (0,1)$ is a parameter of our choice. We may now finish by choosing $c_2$ small enough so that\footnote{We can do so since $\lim_{u\to0^+}I(u) = 1$.} $I(c_2) > 1-\varepsilon/2$ and then $c_1$ large enough so that\footnote{We can do so since $I'(u) = \log u$, whence $I(c_1+c_2)-I(c_1) \geq c_2\log c_1$ by the mean value theorem.} $I(c_1+c_2)-I(c_1)>2$.
\end{proof}

\begin{remark}
    \label{rem:small_h_and_large_t}
    It turns out that Proposition \ref{prop:small_h} admits a proof which at the same time ensures that $T(A) \geq -\varepsilon H(A)\log n$. This means that even the stronger bound
    \[
        \E\bigl[f_A(U)\bigr] \leq \frac{\log H(A) + \frac{T(A)}{H(A)} + 2 - \log\pi}{H(A)},
    \]
    obtained by retaining the $T(A)/H(A)^2$ term in the deduction of Theorem \ref{thm:bound_in_terms_of_h} from Proposition \ref{prop:int_var}, cannot lower the exponent $2$ in Theorem \ref{thm:min_fract_part_new} by any fixed positive amount. We briefly sketch here how such a set $A$ can be constructed. The starting point is to observe that if all elements of $A$ are $y$-smooth, then so are all $d \in D(A)$, so according to Mertens' first theorem \cite[Theorem I.1.7]{tenenbaum}, for all $d \in D(A)$ we have
    \[
        S(d) \leq \sum_{p \leq y}\frac{\log p}{p-1} \leq \log y + O(1).
    \]
    Hence, it follows that
    \[
        -T(A) \leq \sum_{d\in D(A)}\frac{\varphi(d)S(d)}{m_A(d)} \leq H(A)(\log y+O(1)).
    \]
    We now follow the construction from the above proof, by choosing $x = n(\log n)^{I(c_1)+r}$ and $y = (2x)^{1/r}$ for some fixed $r \in \N$ such that $r > 1/\varepsilon$. We also let $k = \lceil c_1\log\log x\rceil$ and $z = x^{1/(2k)}$. Letting $\cD_{u,v,s}$ denote the set of all products of $s$ distinct primes from $(u,v]$, we make the choice
    \[
        B = \bigcup_{d \in \cD_{k^2,z,k-r}}\bigl(d\cdot \cD_{(x/d)^{1/r},(2x/d)^{1/r},r}\bigr).
    \]
    By construction, $B \subseteq (x, 2x]$ and all elements of $B$ are $y$-smooth. Note also that for each $d \in \cD_{k^2,z,k-r}$ we have $d \leq x^{(k-r)/(2k)}$ and hence $(x/d)^{1/r} \geq zx^{1/(2r)} > z$. In particular, each $b \in B$ has $\omega(b) = k$ and the above dilates are pairwise disjoint. Hence, we have
    \[
        |B| = \sum_{d \in \cD_{k^2,z,k-r}}\binom{\pi\bigl((2x/d)^{1/r}\bigr)-\pi\bigl((x/d)^{1/r}\bigr)}{r} \gg_r \frac{x}{(\log x)^r}\sum_{d \in \cD_{k^2,z,k-r}}\frac{1}{d}.
    \]
    To bound the last sum, note that by Mertens' second theorem,
    \[
        \sigma \vcentcolon= \sum_{k^2 < p \leq z}\frac{1}{p} = \log\log z - \log\log(k^2) + O(1) = (1-o(1))\log\log x.
    \]
    Taking $(k-r)$-th powers, expanding out, and discarding the contribution of products containing repeated factors, we obtain
    \[
        \sum_{d \in \cD_{k^2,z,k-r}}\frac{1}{d} \geq \frac{1}{(k-r)!}\Bigl(\sigma^{k-r} - \binom{k-r}{2}\sigma^{k-r-2}\sum_{p>k^2}\frac{1}{p^2}\Bigr) = (1-o(1))\frac{\sigma^{k-r}}{(k-r)!}.
    \]
    Stirling's approximation now gives the estimate $|B| \geq x/(\log x)^{I(c_1)+r-1+o(1)}$, and hence the rest of the proof of Proposition \ref{prop:small_h} goes through; note in particular that all elements of $A$ are $y$-smooth as desired.
\end{remark}

\bibliographystyle{amsplain}
\bibliography{references}

\end{document}